\input amstex
\input amsppt.sty
\magnification=\magstep1
\hsize=32truecc
\vsize=22.2truecm
\baselineskip=16truept
\NoBlackBoxes
\TagsOnRight \pageno=1 \nologo
\def\Z{\Bbb Z}

\def\Q{\Bbb Q}

\def\l{\left}
\def\r{\right}
\def\bg{\bigg}
\def\({\bg(}
\def\[{\bg\lfloor}
\def\){\bg)}
\def\]{\bg\rfloor}
\def\t{\text}
\def\f{\frac}

\def\se {\subseteq}

\def\bi{\binom}
\def\eq{\equiv}

\def\ls{\leqslant}
\def\gs{\geqslant}
\def\mo{\roman{mod}}

\def\ve{\varepsilon}
\def\al{\alpha}

\def\Proof{\noindent{\it Proof}}
\def\Def{\medskip\noindent{\it Definition}}

\def\Remark{\medskip\noindent{\it  Remark}}

\def\Ack{\medskip\noindent {\bf Acknowledgments}}
\hbox {Int. J. Number Theory 12 (2016), no.8, 2323--2342.}
\bigskip
\topmatter
\title On a pair of zeta functions \endtitle
\author Zhi-Wei Sun\endauthor
\leftheadtext{Zhi-Wei Sun}
 \rightheadtext{On a pair of zeta functions}
\affil Department of Mathematics, Nanjing University\\
 Nanjing 210093, People's Republic of China
  \\  zwsun\@nju.edu.cn
  \\ {\tt http://math.nju.edu.cn/$\sim$zwsun}
\endaffil
\abstract Let $m$ be a positive integer, and define
$$\zeta_m(s)=\sum_{n=1}^\infty\f{(-e^{2\pi i/m})^{\omega(n)}}{n^s}\ \ \t{and}
\ \ \zeta^*_m(s)=\sum_{n=1}^\infty\f{(-e^{2\pi i/m})^{\Omega(n)}}{n^s},$$
for $\Re(s)>1$, where $\omega(n)$ denotes the number of distinct prime factors of $n$, and $\Omega(n)$
represents the total number of prime factors of $n$ (counted with multiplicity).
In this paper we study these two zeta functions and related arithmetical functions.
We show that
$$\sum^\infty\Sb n=1\\ n\ \t{is squarefree}\endSb\f{(-e^{2\pi i/m})^{\omega(n)}}n=0\quad\t{if}\ m>4,$$
which is similar to the known identity
$\sum_{n=1}^\infty\mu(n)/n=0$ equivalent to the Prime Number
Theorem. For $m>4$, we prove that
$$\zeta_m(1):=\sum_{n=1}^\infty\f{(-e^{2\pi i/m})^{\omega(n)}}n=0
\ \ \t{and}\ \ \zeta^*_m(1):=\sum_{n=1}^\infty\f{(-e^{2\pi
i/m})^{\Omega(n)}}n=0.$$
We also raise a hypothesis on the parities of $\Omega(n)-n$ which implies the Riemann Hypothesis.

\endabstract
\thanks 2010 {\it Mathematics Subject Classification}.\,Primary 11M99;
Secondary 11A25, 11N37.
\newline\indent Supported by the National Natural Science
Foundation of China (grant 11571162) .
\endthanks
\endtopmatter
\document

\heading{1. Introduction}\endheading

The Riemann zeta function $\zeta(s)$, defined by
$$\zeta(s)=\sum_{n=1}^\infty\f1{n^s}\ \ \quad\t{for}\ \Re(s)>1,$$
plays a very important role in number theory. As Euler observed,
$$\zeta(s)=\prod_p\l(1-\f1{p^s}\r)^{-1}\ \ \t{for}\ \Re(s)>1.$$
(In such a product we always let $p$ run over all primes.)
It is well-known that $\zeta(s)$ for $\Re(s)>1$ can be continued analytically to a complex function which is holomorphic everywhere except for a simple
pole at $s=1$ with residue 1. The famous Riemann Hypothesis asserts that
if $0\ls\Re(s)\ls 1$ and $\zeta(s)=0$ then $\Re(s)=1/2$.
The Prime Number Theorem $\pi(x)\sim x/\log x$ (as $x\to+\infty$)
is actually equivalent to $\zeta(1+it)\not=0$ for any nonzero real number $t$.
(See, e.g., R. Crandall and C. Pomerance [CP, pp.\,33-37].)

The M\"obius fucntion $\mu$ defined on $\Z^+=\{1,2,3,\ldots\}$ is given by
$$\mu(n)=\cases1&\t{if}\ n=1,
\\(-1)^k&\t{if}\ n\ \t{is a product of}\ k\ \t{distinct primes},
\\0&\t{if}\ p^2\mid n\ \t{for some prime}\ p.
\endcases$$
It is well known that
$$\zeta(s)\sum_{n=1}^\infty\f{\mu(n)}{n^s}=1\ \ \t{for}\ \Re(s)>1.$$
Also, either of $\sum_{n=1}^\infty\mu(n)/n=0$ and $\sum_{n\ls x}\mu(n)=o(x)$ is equivalent to the Prime Number Theorem.
(Cf. T. M. Apostol [Ap, \S3.9 and \S4.1].)

The reader may consult [Ap] and [IR, pp.\,18-21] for the basic knowledge of arithmetical functions
and the theory of Dirichlet's convolution and Dirichlet series.

If $n\in\Z^+$ is squarefree, then $\mu(n)=(-1)^{\Omega(n)}$ depends on $\Omega(n)$ modulo 2, where $\Omega(n)$ denotes the number
of all prime factors of $n$ (counted with multiplicity). For the Liouville function $\lambda(n)=(-1)^{\Omega(n)}$, it is known that
$$\sum_{d\mid n}\lambda(d)=\cases1&\t{if}\ n\ \t{is a square},\\0&\t{otherwise}.\endcases$$
(See, e.g., [Ap, p.\,38].)
J. van de Lune and R. E. Dressler [LD] showed that $\sum_{n=1}^\infty(-1)^{\omega(n)}/n=0$,
where $\omega(n)$ denotes the number of distinct prime factors of $n$.

Now we give natural extensions of the functions $\mu(n)$, $\lambda(n)$ and $\zeta(s)$.

\Def\ 1.1. Let $m$ be any positive integer. For $n\in\Z^+$ we set
$$\mu_m(n)=\cases (-e^{2\pi i/m})^{\omega(n)}&\t{if}\ n\ \t{is squarefree},
\\0&\t{otherwise},\endcases\tag1.1$$
$$\nu_m(n)=(-e^{2\pi i/m})^{\omega(n)}\ \ \t{and}\ \ \nu_m^*(n)=(-e^{2\pi i/m})^{\Omega(n)}.\tag1.2$$
For $\Re(s)>1$ we define
$$\zeta_m(s)=\sum_{n=1}^\infty\f{\nu_m(n)}{n^s}=\prod_p\l(1-\f{e^{2\pi i/m}}{p^s-1}\r)\tag1.3$$
and
$$\zeta_m^*(s)=\sum_{n=1}^\infty\f{\nu_m^*(n)}{n^s}=\prod_p\l(1+\f{e^{2\pi i/m}}{p^s}\r)^{-1}.\tag1.4$$
\medskip

As $\nu_m^*$ is completely multiplicative, the second identity in (1.4) is easy and in fact known.
Since $\nu_m$ is multiplicative, if $\Re(s)>1$ then
$$\sum_{n=1}^\infty\f{\nu_m(n)}{n^s}=\prod_p\sum_{k=0}^\infty\f{\nu_m(p^k)}{p^{ks}}=\prod_p\bg(1-e^{2\pi i/m}\sum_{k=1}^\infty\f1{p^{ks}}\bg)$$
and hence the second equality in (1.3) does hold.

As $\mu_1=\mu$, we call $\mu_m$ the generalized M\"obius function of order $m$.
Note that $\zeta_2(s)=\zeta_2^*(s)=\zeta(s)$.
Also, $\nu_1^*(n)=(-1)^{\Omega(n)}$ is the Liouville function $\lambda(n)$, and
$$\zeta_1^*(s)=\sum_{n=1}^\infty\f{\lambda(n)}{n^s}=\f{\zeta(2s)}{\zeta(s)}=\prod_p\l(1+\f1{p^s}\r)^{-1}\ \ \t{for}\ \Re(s)>1.$$
(Cf. [Ap, pp.\,229-230].)
If we replace $-e^{2\pi i/m}$ in the definition of $\zeta_m^*(s)$ by $e^{2\pi i/m}$, the resulting function
was shown to have an infinitely many valued analytic continuation into the half plane $\Re(s)>1/2$ by T. Kubota and M. Yoshida [KY].
(See also [A] and [CD].) It seems that the zeta function $\zeta_m(s)$ introduced here has not been studied before.

Our first theorem is a basic result.

\proclaim{Theorem 1.1} Let $m$ be any positive integer.

{\rm (i)} The function $\mu_m^*(n)=\mu_m(n)\lambda(n)$ is the inverse of $\nu^*_m(n)$ with respect to the Dirichlet convolution, and hence
$$\zeta^*_m(s)\sum_{n=1}^\infty\f{\mu^*_m(n)}{n^s}=1\qquad \t{for}\ \Re(s)>1.\tag1.5$$
For $\Re(s)>1$ we also have
$$\zeta_m(s)\sum_{n=1}^\infty\f{(1+e^{2\pi i/m})^{\Omega(n)}}{n^s}=\zeta(s).\tag1.6$$

{\rm (ii)} If $m>4$, then
$$\prod_p\l(1+\f{e^{2\pi i/m}}p\r)^{-1}=0.\tag1.7$$
On the other hand,
$$\prod_p\l(1+\f{e^{2\pi i/3}}p\r)=0\ \ \t{and}\ \ \lim_{x\to\infty}\bigg|\prod_{p\ls x}\l(1+\f{e^{2\pi i/4}}p\r)\bigg|=\f{\sqrt{15}}{\pi}.\tag1.8$$
\endproclaim
\Remark\ 1.1.  If $\Re(s)>1$, then both $\zeta_m^*(s)$ and $\zeta_m(s)$ are nonzero by (1.5) and (1.6).
\medskip

Our second theorem is a general result.
\proclaim{Theorem 1.2} Let $z$ be a complex number with $\Re(z)<1$. For $x\gs2$ we have
$$\sum_{n\ls x}\f{z^{\omega(n)}}n=\Cal F(z)(\log x)^z+c(z)+O((\log x)^{z-1})\tag1.9$$
and
$$\sum\Sb n\ls x\\ n\ \t{is squarefree}\endSb\f{z^{\omega(n)}}n=\Cal G(z)(\log x)^z+c_*(z)+O((\log x)^{z-1}),\tag1.10$$
where $c(z)$ and $c_*(z)$ are constants only depending on $z$, and
$$\align \Cal F(z)=&\f1{\Gamma(1+z)}\prod_p\l(1+\f z{p-1}\r)\l(1-\f1p\r)^z,
\\\Cal G(z)=&\f1{\Gamma(1+z)}\prod_p\l(1+\f z{p}\r)\l(1-\f1p\r)^z.
\endalign$$
If $|z|<2$, then for $x\gs 2$ we have
$$\sum_{n\ls x}\f{z^{\Omega(n)}}n=\Cal H(z)(\log x)^z+C(z)+O((\log x)^{z-1}),\tag1.11$$
where $C(z)$ is a constant only depending on $z$, and
$$\Cal H(z)=\f1{\Gamma(1+z)}\prod_p\l(1-\f zp\r)^{-1}\l(1-\f1p\r)^z.$$
\endproclaim

Theorem 1.2 obviously has the following consequence.
\proclaim{Corollary 1.1} For any complex number $z$ with $\Re(z)<0$, we have
$$\sum_{n=1}^\infty\f{z^{\omega(n)}}n=c(z)\ \ \t{and}\ \ \sum^\infty\Sb n=1\\ n\ \t{is squarefree}\endSb\f{z^{\omega(n)}}n=c_*(z).\tag1.12$$
If $|z|<2$ and $\Re(z)<0$, then
$$\sum_{n=1}^\infty\f{z^{\Omega(n)}}n=C(z).\tag1.13$$
\endproclaim

\proclaim{Theorem 1.3} We have
$$\sum_{n=1}^\infty\f{\mu_5(n)}n=\sum_{n=1}^\infty\f{\mu_6(n)}n=\cdots=0.\tag1.14$$
Moreover, for any positive integer $m\not=2$ we have
$$(\log x)^{e^{2\pi i/m}}\sum_{n\ls x}\f{\mu_m(n)}n=\Cal G(-e^{2\pi i/m})+O\l(\f1{\log x}\r)\ \ (x\gs2),\tag1.15$$
where $\Cal G(z)$ is defined as in Theorem 1.2.
\endproclaim
\Remark\ 1.2. It is known that
$$\sum_{n\ls x}\f{\mu_2(n)}n=\sum_{n\ls x}\f{|\mu(n)|}n=\f 6{\pi^2}\log x+c+O\l(\f1{\sqrt x}\r)\ (x\gs 2),$$
where $c=1.04389\ldots$ (see, e.g., [BS, Lemma 14]).
(1.15) with $m=4$ implies that
$$\lim_{x\to\infty}\bigg|\sum_{n\ls x}\f{\mu_4(n)}n\bigg|=|\Cal G(-i)|.$$
After reading the first version of this paper, D. Broadhurst simplified $|\Cal G(-i)|$ as $\sqrt{15(\sinh\pi)/\pi^3}.$

\proclaim{Theorem 1.4} Let
$$V_m(x)=\sum_{n\ls x}\f{\nu_m(n)}n\ \ \t{and}\ \ V^*_m(x)=\sum_{n\ls x}\f{\nu^*_m(n)}n$$
 for $m\in\Z^+$ and $x\gs 2$. Then
$$\aligned V_3(x)=&\Cal F(-e^{2\pi i/3})(\log x)^{(1-i\sqrt3)/2}+c_3+O\l(\f1{\sqrt{\log x}}\r),
\\V^*_3(x)= &\Cal H(-e^{2\pi i/3})(\log x)^{(1-i\sqrt{3})/2}+C_3+O\l(\f1{\sqrt{\log x}}\r),\endaligned\tag1.16$$
and
$$\aligned V_4(x)=&\Cal F(-i)(\log x)^{-i}+c_4+O\l(\f1{\log x}\r),
\\V^*_4(x)= &\Cal H(-i)(\log x)^{-i}+C_4+O\l(\f1{\log x}\r),\endaligned\tag1.17$$
where $c_3,C_3,c_4,C_4$ are suitable constants.
Also, for $m=5,6,\ldots$ we have $V_m(x)=o(1)$ and $V_m^*(x)=o(1)$, i.e.,
$$\zeta_m(1):=\sum_{n=1}^\infty\f{\nu_m(n)}n=0\ \ \t{and}\ \ \zeta^*_m(1):=\sum_{n=1}^\infty\f{\nu^*_m(n)}n=0.\tag1.18$$
Moreover, for $m=1,5,6,\ldots$ we have
$$V_m(x)(\log x)^{e^{2\pi i/m}}=\Cal F(-e^{2\pi i/m})+O\l(\f1{\log x}\r)\tag1.19$$
and
$$V^*_m(x)(\log x)^{e^{2\pi i/m}}=\Cal H(-e^{2\pi i/m})+O\l(\f1{\log x}\r).\tag1.20$$
\endproclaim
\Remark\ 1.3. It seems that $c_3$ and $C_3$ are nonzero but $c_4=0$ (and probably also $C_4=0$).  Broadhurst simplified $|\Cal H(-i)|$ as $\sqrt{(\sinh\pi)\pi/15}$.
\medskip

Theorem 1.1 is not difficult. Our proofs of Theorems 1.2-1.4 depend heavily on some results of A. Selberg [S] (see also H. Delange [D]
and Theorem 7.18 of [MV, p.\,231]) and the partial summation method via Abel's identity (see, [Ap, p.\,77]).

Motivated by Theorem 1.4 we pose the following conjecture for further research.

\proclaim{Conjecture 1.1} Both $V_1(x)=\sum_{n\ls x}(-1)^{\omega(n)}/n$
and $V_1^*(x)=\sum_{n\ls x}(-1)^{\Omega(n)}/n$ are $O(x^{\ve-1/2})$ for any $\ve>0$.
Also, $|\sum_{n\ls x}(-2)^{\Omega(n)}|<x$ for all $x\gs 3078$.
\endproclaim
\Remark\ 1.4. It seems that $V_1(x)$ might be $O(\sqrt{(\log x)/x})$ or even $O(1/\sqrt x)$.
The asymptotic behavior of  $\sum_{n\ls x}2^{\Omega(n)}$ was investigated by E. Grosswald [G].

\medskip
In 1958 C. B. Haselgrove [H] disproved P\'olya's conjecture that $\sum_{n\ls x}\lambda(n)\ls0$ for all $x\gs2$; he also showed that
Tur\'an's conjecture $\sum_{n\ls x}\lambda(n)/n$ $>0$ for $x\gs1$, is also false.
It is known that the least integer $x>1$ with $\sum_{n\ls x}\lambda(n)>0$
is $906150257 < 10^9$ (cf. [L] and [BFM]).
Along this line we propose the following new hypothesis.

\proclaim{Hypothesis 1.1} {\rm (i)} For any $x\gs 5$, we have
$$S(x):=\sum_{n\ls x}(-1)^{n-\Omega(n)}>0,\tag1.21$$
i.e.,
$$|\{n\ls x:\ \Omega(n)\eq n\ (\mo\ 2)\}|>|\{n\ls x:\ \Omega(n)\not\eq n\ (\mo\ 2)\}|.$$
Moreover,
$$S(x)>\sqrt x \ \ \t{for all}\ x\gs 325,\ \t{and}\ \ S(x)<2.3\sqrt x\ \ \t{for all}\ x\gs1.$$

{\rm (ii)} For any $x\gs 1$ we have
$$T(x):=\sum_{n\ls x}\f{(-1)^{n-\Omega(n)}}n<0.\tag1.22$$
Moreover,
$$T(x)\sqrt x<-1\ \ \t{for all}\ x\gs 2,\ \ \t{and}\ \
T(x)\sqrt x>-2.3\ \ \t{for all}\ x\gs 3.$$
\endproclaim
\Remark\ 1.5. We have verified parts (i) and (ii) of the hypothesis for $x$ up to $10^{11}$ and $2\times10^9$ respectively.
Below are values of $S(x)$ for some particular $x$:
$$\gather S(10^2)=14,\  S(10^3)=54,\ S(10^4)=186,\ S(10^5)=464,\ S(10^6)=1302,
\\\ S(10^7)=5426,\ S(10^8)=19100,\  S(10^9)=62824,\ S(10^{10})=172250,
\\ S(2\cdot10^{10})=252292,\  S(3\cdot10^{10})=292154,\ S(4\cdot10^{10})=263326,
\\ S(5\cdot10^{10})=360470,\ S(6\cdot10^{10})=363152,\ S(7\cdot10^{10})=406260,
\\ S(8\cdot10^{10})=559558,\ S(9\cdot10^{10})=491100,\ S(10^{11})=457588.
\endgather$$

{\it Example}\ 1.1. For $x_1=17593752$ and $x_2=123579784$, we have
$S(x_1)=9574$ and $S(x_2)=11630$. Via a computer we find that
$$\max_{1\ls x\ls10^{11}}\f{S(x)}{\sqrt x}=\f{S(x_1)}{\sqrt{x_1}}\approx 2.28252$$
and
$$\min_{324<x\ls10^{11}}\f{S(x)}{\sqrt x}=\ \f{S(x_2)}{\sqrt{x_2}}\approx 1.04618.$$
\medskip

We are unable to prove or disprove Hypothesis 1.1, but we can show the following relatively easy result.

\proclaim{Theorem 1.5} {\rm (i)} We have
$$S(x)=o(x)\ \ \t{and}\ \ \sum_{n=1}^\infty\f{(-1)^{n-\Omega(n)}}n=0.\tag1.23$$

{\rm (ii)} If $S(x)>0$ for all $x\gs 5$, or $T(x)<0$ for all $x\gs1$, then the Riemann Hypothesis holds.
\endproclaim

Note that
$$S(x)>0\iff|\{n\ls x:\ 2\mid n-\Omega(n)\}|>\f x2.$$
In view of Hypothesis 1.1, it is natural to ask whether
$$|\{n\ls x:\ m\mid n-\Omega(n)\}|>\f xm \ \t{for sufficiently large}\ x.$$
For $m=3,4,\ldots,18,20$ we have the following conjecture based on our computation.

\proclaim{Conjecture 1.2} We have
$$|\{n\ls x: 4\mid n-\Omega(n)\}| < \f x4\quad\t{for any}\  x\gs s(4),$$
and for $m=3,5,6,\cdots,18,20$ we have
 $$ |\{n\ls x: m\mid n-\Omega(n)\}| > \f xm\quad\t{for all}\  x\gs s(m),$$
where
 $$\align&s(3)=62,\ s(4)=1793193,\ s(5)=187,\ s(6)=14,\ s(7)=6044,\ s(8)=73,
 \\&s(9)=65,\ s(10)=61,\ s(11)=4040389,\ s(12)=14,\ s(13)=6943303,
 \\&s(14)=4174,\ s(15)=77,\ s(16)=99,\ s(17)=50147927,\ s(18)=73,\ s(20)=61.
 \endalign$$
\endproclaim
\Remark\ 1.7.  The case $m=19$ seems much more sophisticated. Perhaps the sign of $|\{n\ls x: 19|(n-\Omega(n))\}|-x/19$ changes infinitely often.
\medskip

As there is an extended Riemann Hypothesis for algebraic number fields, we propose the following extension of Hypothesis 1.1
based on our computation.

\proclaim{Hypothesis 1.2 (Extended Hypothesis)} Let $K$ be any algebraic number field.
Then we have
$$S_K(x):=\sum_{N(A)\ls x}(-1)^{N(A)-\Omega(A)}>0\quad\t{for all sufficiently large}\ x,$$
where $A$ runs over all nonzero integral ideals in $K$ whose norm (with respect to the field extension $K/\Q$) are not greater than $x$, and $\Omega(A)$ denotes the total number of prime ideals
in the factorization of $A$ as a product of prime ideals (counted with multiplicity).
In particular, for $K=\Q(i)$ we have $S_K(x)>0$ for all $x\gs 9$, and for $K=\Q(\sqrt{-2})$ we have $S_K(x)>0$ for all $x\gs 132.$
\endproclaim

Now we give one more conjecture based on our computation.
\proclaim{Conjecture 1.3} For an integer $d\eq0,1\ (\mo\ 4)$ define
$$S_d(x)=\sum_{n\ls x}(-1)^{n-\Omega(n)}\l(\f dn\r),$$
where $(\f dn)$ denotes the Kronecker symbol. Then
$$S_{-4}(x)<0,\ S_{-7}(x)<0,\ S_{-8}(x)<0$$
for all $x\gs 1$, and
$$S_5(x)>0\ \t{for}\ x\gs 11,\ \ S_{-3}(x)>0\ \t{for}\ x\gs 406759,\ S_{-11}(x)>0\ \t{for}\ x\gs771862,$$
and
$$S_{24}(x)<0\ \t{for}\ x\gs90601,\ \ \t{and}\ \ S_{28}(x)<0\ \t{for}\ x\gs 629819.$$
\endproclaim

We will show Theorems 1.1 and 1.2 in the next section, and prove Theorems 1.3-1.5 in Sections 3-5 respectively.

\heading{2. Proofs of Theorems 1.1 and 1.2}\endheading

\medskip\noindent
{\it Proof of Theorem 1.1}. (i) Clearly $\mu^*_m(1)\nu_m^*(1)=1\cdot1=1$. Let $N$ be any integer greater than one, and let $n$ be the product of all distinct prime factors of $N$.
Then
$$\align\sum_{d\mid N}\mu^*_m(d)\nu^*_m\l(\f Nd\r)=&\sum_{d\mid n}e^{2\pi i\Omega(d)/m}(-e^{2\pi i/m})^{\Omega(n/d)+\Omega(N/n)}
\\=&(-1)^{\Omega(N/n)}e^{2\pi i\Omega(N)/m}\sum_{d\mid n}\mu\l(\f nd\r)=0.
\endalign$$
Therefore $\mu^*_m$ is the inverse of $\nu^*_m$ with respect to the Dirichlet convolution $*$.

Let $s=\sigma+it$ be a complex number with $\Re(s)=\sigma>1$. Since
$$\max\l\{\l|\f{\mu^*_m(n)}{n^s}\r|,\l|\f{\nu^*_m(n)}{n^s}\r|\r\}\ls\l|\f{1}{n^{\sigma+it}}\r|=\l|\f{e^{-it\log n}}{n^\sigma}\r|=\f1{n^{\sigma}}$$
for any $n\in\Z^+$, both $\sum_{n=1}^\infty\mu^*_m(n)/n^s$ and $\sum_{n=1}^\infty\nu^*_m(n)/n^s$
converge absolutely. Therefore
$$\zeta_m^*(s)\sum_{n=1}^\infty\f{\mu^*_m(n)}{n^s}=\sum_{n=1}^\infty\f{\mu^*_m(n)}{n^s}\sum_{n=1}^\infty\f{\nu^*_m(n)}{n^s}=\sum_{n=1}^\infty\f{\mu_m^**\nu^*_m(n)}{n^s}=1.$$

Now we prove (1.6). Since $|p^s|=p^{\sigma}>p\gs|1+e^{2\pi i/m}|$ for any prime $p$, we have
$$\prod_p\l(1-\f{1+e^{2\pi i/m}}{p^s}\r)^{-1}=\prod_p\sum_{k=0}^\infty\f{(1+e^{2\pi i/m})^k}{p^{ks}}
=\sum_{n=0}^\infty\f{(1+e^{2\pi i/m})^{\Omega(n)}}{n^s}.$$
Note that
$$\align\zeta_m(s)=&\prod_p\f{p^s-1-e^{2\pi i/m}}{p^s-1}=\prod_p\f{1-(1+e^{2\pi i/m})/p^s}{1-1/p^s}
\\=&\zeta(s)\prod_p\l(1-\f{1+e^{2\pi i/m}}{p^s}\r).\endalign$$
So (1.6) does hold.

(ii) Now assume that $m>4$. Then $2\pi/m<\pi/2$ and $0<\cos(2\pi/m)<1$. For any prime $p$ we have
$$\l|1+\f{e^{2\pi i/m}}p\r|=\l|\l(1+\f{\cos(2\pi/m)}p\r)+i\f{\sin(2\pi /m)}p\r|\gs 1+\f{\cos(2\pi/m)}p.$$
Therefore
$$\bigg|\prod_{p\ls x}\l(1+\f{e^{2\pi i/m}}p\r)\bigg|\gs\prod_{p\ls x}\l(1+\f{\cos(2\pi/m)}p\r)\gs 1+\cos\f{2\pi}m\sum_{p\ls x}\f1p,$$
and hence (1.7) holds since $\sum_p1/p$ diverges (cf. [IR, p.\,21]).

Finally we prove the first identity in (1.8). For any prime $p$, we have
$$\l|1+\f{e^{2\pi i/3}}p\r|^2=1+2\f{\cos{2\pi/3}}p+\f1{p^2}=1-\f1p+\f1{p^2}=\f{1+p^{-3}}{1+p^{-1}}.$$
Thus
$$\align\bg|\prod_{p\ls x}\l(1+\f{e^{2\pi i/3}}p\r)\bg|^2=&\prod_{p\ls x}\l(1+\f1{p^3}\r)\cdot\prod_{p\ls x}\l(1+\f1p\r)^{-1}
\\\ls&\prod_p\l(1+\f1{p^3}\r)\cdot\bg(1+\sum_{p\ls x}\f1p\bg)^{-1}.
\endalign$$
Since $\sum_p1/p$ diverges while $\sum_p 1/p^3$ converges, the first equality in (1.8) follows.

The second equality in (1.8) is easy. In fact, as $x\to\infty$,
$$\bigg|\prod_{p\ls x}\l(1+\f{e^{2\pi i/4}}p\r)\bigg|^2=\prod_{p\ls x}\l|1+\f ip\r|^2$$
has the limit
$$\prod_p\l(1+\f1{p^2}\r)=\f{\prod_p(1-1/p^2)^{-1}}{\prod_p(1-1/p^4)^{-1}}=\f{\zeta(2)}{\zeta(4)}=\f{\pi^2/6}{\pi^4/90}=\f{15}{\pi^2}.$$

In view of the above, we have completed the proof of Theorem 1.1. \qed
\medskip

To prove Theorem 1.2, we need two lemmas.

\proclaim{Lemma 2.1 {\rm (Selberg [S])}} Let $z$ be a complex number. For $x\gs2$ we have
$$\sum_{n\ls x}z^{\omega(n)}=F(z)x(\log x)^{z-1}+O\l(x(\log x)^{\Re(z)-2}\r)\tag2.1$$
and
$$\sum\Sb n\ls x\\n\ \t{is squarefree}\endSb z^{\omega(n)}=G(z)x(\log x)^{z-1}+O\l(x(\log x)^{\Re(z)-2}\r),\tag2.2$$
where
$$F(z)=\f1{\Gamma(z)}\prod_p\l(1+\f z{p-1}\r)\l(1-\f1p\r)^z$$
and
$$G(z)=\f1{\Gamma(z)}\prod_p\l(1+\f z{p}\r)\l(1-\f1p\r)^z.$$
When $|z|<2$, for $x\gs2$ we also have
$$\sum_{n\ls x}z^{\Omega(n)}=H(z)x(\log x)^{z-1}+O\l(x(\log x)^{\Re(z)-2}\r),\tag2.3$$
where
$$H(z)=\f1{\Gamma(z)}\prod_p\l(1-\f zp\r)^{-1}\l(1-\f1p\r)^z.$$
\endproclaim

\proclaim{Lemma 2.2} Let $a(1),a(2),\ldots$ be a sequence of complex numbers. Suppose that
$$\sum_{n\ls x}a(n)=cx(\log x)^{z-1}+O(x(\log x)^{\Re(z)-2})\ \ (x\gs2),\tag2.4$$
where $c$ and $z$ are (absolute) complex numbers with $z\not=0$ and $\Re(z)\not=1$. Then, for $x,y\gs2$ we have
$$\aligned&\sum_{n\ls x}\f{a(n)}n-\f{c}{z}(\log x)^z-\(\sum_{n\ls y}\f{a(n)}n-\f{c}{z}(\log y)^z\)
\\&\quad=O((\log x)^{z-1})+O((\log y)^{z-1}).
\endaligned\tag2.5$$
Thus, if $\Re(z)<1$ then
$$\sum_{n\ls x}\f{a(n)}n=\f{c}{z}(\log x)^z+c_z+O((\log x)^{\Re(z)-1})\ \ (x\gs2),\tag2.6$$
where $c_z$ is a suitable constant.
\endproclaim
\Proof.  Let $A(t)=\sum_{n\ls t}a(n)$ for $t\gs2$. By the partial summation formula,
$$\align\sum_{n\ls x}\f{a(n)}n-\sum_{n\ls y}\f{a(n)}n=&\f{A(x)}x-\f{A(y)}{y}-\int_{y}^x A(t)(t^{-1})'dt
\\=&\f{A(x)}x-\f{A(y)}{y}+\int_{y}^x \f{A(t)}{t^2}dt.
\endalign$$
Note that
$$\f{A(t)}t=c(\log t)^{z-1}+O((\log t)^{\Re(z)-2})\quad\t{for}\ t\gs 2.$$
Clearly
$$\int_{y}^x\f{(\log t)^{z-1}}tdt=\f{(\log t)^{z}}{z}\bg|_{t=y}^x=\f{(\log x)^{z}-(\log y)^{z}}{z}$$
and
$$\int_{y}^x\f{(\log t)^{\Re(z)-2}}tdt=\f{(\log t)^{\Re(z)-1}}{\Re(z)-1}\bg|_{t=y}^x=\f{(\log x)^{\Re(z)-1}-(\log y)^{\Re(z)-1}}{\Re(z)-1}.$$
So the desired (2.5) follows from the above.

Now assume that $\Re(z)<1$. For any $\ve>0$ we can find a positive integer $N$ such that for $x,y\gs N$
the absolute value of the right-hand side of (2.5) is smaller than $\ve$. Therefore, in view of (2.5) and Cauchy's convergence criterion,
 $\sum_{n\ls x}a(n)/n-c(\log x)^z/z$ has a finite limit $c_z$
as $x\to\infty$. Letting $y\to\infty$ in (2.5) we immediately obtain (2.6). This ends the proof. \qed

\medskip
\noindent{\it Proof of Theorem 1.2}. When $z=0$, (1.9)-(1.11) obviously hold with $c(0)=c_*(0)=C(0)=0$.

Now assume $z\not=0$. As $\Gamma(1+z)=z\Gamma(z)$, we see that
$$\Cal F(z)=\f{F(z)}z,\ \Cal G(z)=\f{G(z)}z,\ \t{and}\ \Cal H(z)=\f{H(z)}z,$$
where the functions $F$, $G$ and $H$ are given in Lemma 2.1.
Combining Lemmas 2.1 and 2.2 we immediately get the desired (1.9)-(1.11). \qed

\heading{3. Proof of Theorem 1.3}\endheading

We first present two lemmas.

\proclaim{Lemma 3.1} Let $m\in\Z^+$ and $x\gs1$. Then we have
$$\sum_{n\ls x}\mu_m(n)\l\lfloor\f xn\r\rfloor=\sum_{n\ls x}(1-e^{2\pi i/m})^{\omega(n)}.\tag3.1$$
\endproclaim
\Proof. We first claim that
$$\sum_{d\mid n}\mu_m(d)=(1-e^{2\pi i/m})^{\omega(n)}\tag 3.2$$
for any $n\in\Z^+$. Clearly (3.2) holds for $n=1$. If $n=p_1^{a_1}\cdots p_k^{a_k}$ with $p_1,\ldots,p_k$ distinct primes and $a_1,\ldots,a_k\in\Z^+$,
then
$$\sum_{d\mid n}\mu_m(d)=\sum_{I\se\{1,\ldots,k\}}\mu_m\bg(\prod_{i\in I}p_i\bg)=\sum_{r=0}^k\bi kr(-e^{2\pi i/m})^r=(1-e^{2\pi i/m})^{\omega(n)}.$$

Observe that
$$\sum_{d\ls x}\mu_m(d)\l\lfloor\f xd\r\rfloor=\sum_{d\ls x}\mu_m(d)\sum_{q\ls x/d}1=\sum_{dq\ls x}\mu_m(d)=\sum_{n\ls x}\sum_{d\mid n}\mu_m(d).$$
Combining this with (3.2) we immediately obtain (3.1). \qed

\proclaim{Lemma 3.2} Let $m\in\Z^+$, $m\not=2$, and $x\gs2$. Then we have
$$\sum_{n\ls x}\mu_m(n)\l\{\f xn\r\}=o(x),\ \sum_{n\ls x}\nu_m(n)\l\{\f xn\r\}=o(x),\ \sum_{n\ls x}\nu^*_m(n)\l\{\f xn\r\}=o(x),\tag3.3$$
where $\{\al\}$ denotes the fractional part of a real number $\al$.
\endproclaim
\Proof. As $m\not=2$, $\Re(e^{2\pi i/m})=\cos\f{2\pi}m\not=-1$. Applying (2.1)-(2.3) we obtain
$$\align\sum_{n\ls x}\mu_m(x)=&xG(-e^{2\pi i/m})(\log x)^{-e^{2\pi i/m}-1}+O\l(x(\log x)^{-\cos(2\pi/m)-2}\r)=o(x),
\\\sum_{n\ls x}\nu_m(x)=&xF(-e^{2\pi i/m})(\log x)^{-e^{2\pi i/m}-1}+O\l(x(\log x)^{-\cos(2\pi/m)-2}\r)=o(x),
\\\sum_{n\ls x}\nu^*_m(x)=&xH(-e^{2\pi i/m})(\log x)^{-e^{2\pi i/m}-1}+O\l(x(\log x)^{-\cos(2\pi/m)-2}\r)=o(x).
\endalign$$
(Note that $F(-1)=G(-1)=H(-1)=0$.)

Let $w$ be any of the three functions $\mu_m$, $\nu_m$ and $\nu^*_m$. By the above, $W(x)=\sum_{n\ls x}w(n)=o(x)$. We want to show that
$$\Delta(x):=\sum_{n\ls x}w(n)\l\{\f xn\r\}=o(x).$$
Clearly
$$r(u):=\sup_{t\gs u}\f{|W(t)|}t\ls1\quad\t{for}\ u\gs1.$$
Also, $r(u)\to0$ as $u\to\infty$.

 Let $0<\ve<1$. Then
 $$\align|\Delta(x)|\ls&\bigg|\sum_{n\ls\ve x}w(n)\l\{\f xn\r\}\bigg|+\bigg|\sum_{\ve x<n\ls x}w(n)\l\{\f xn\r\}\bigg|
\\\ls &\ve x+\bigg|\sum_{\ve x<n\ls x}(W(n)-W(n-1))\l\{\f xn\r\}\bigg|
\\\ls&\ve x+\bigg|\sum_{\ve x<n<\lfloor x\rfloor}W(n)\l(\l\{\f xn\r\}-\l\{\f x{n+1}\r\}\r)\bigg|
\\&+\l|W(\lfloor x\rfloor)\l\{\f x{\lfloor x\rfloor}\r\}-W(\lfloor\ve x\rfloor)\l\{\f x{\lfloor\ve x\rfloor+1}\r\}\r|.
\endalign$$
Note that
$$\l|W(\lfloor x\rfloor)\l\{\f x{\lfloor x\rfloor}\r\}\r|=|W(\lfloor x\rfloor)|\f{\{x\}}{\lfloor x\rfloor}\ls1$$
and
$$\l|W(\lfloor\ve x\rfloor)\l\{\f x{\lfloor\ve x\rfloor+1}\r\}\r|\ls|W(\lfloor\ve x\rfloor)|\ls\lfloor\ve x\rfloor\ls \ve x.$$
Therefore
$$\align|\Delta(x)|\ls& 1+2\ve x+\sum_{\ve x<n<\lfloor x\rfloor}\f{|W(n)|}nx\l|\l\{\f xn\r\}-\l\{\f x{n+1}\r\}\r|
\\\ls&1+2\ve x+xr(\ve x)\sum_{\ve x<n<\lfloor x\rfloor}\bigg |\f xn-\f x{n+1}-\l(\l\lfloor \f x{n}\r\rfloor-\l\lfloor\f x{n+1}\r\rfloor\r)\bigg|
\\\ls&1+2\ve x+xr(\ve x)\sum_{\ve x<n<\lfloor x\rfloor}\l(\l(\f xn-\f x{n+1}\r)+\l(\l\lfloor \f x{n}\r\rfloor-\l\lfloor\f x{n+1}\r\rfloor\r)\r)
\\\ls&1+2\ve x+xr(\ve x)\l(2\f x{\lfloor\ve x\rfloor+1}-\f x{\lfloor x\rfloor}-\l\lfloor\f x{\lfloor x\rfloor}\r\rfloor\r)
\endalign$$
and hence
$$\f{|\Delta(x)|}x\ls \f1x+2\ve+\f 2{\ve}r(\ve x).$$
It follows that
$$\limsup_{x\to\infty}\f{|\Delta(x)|}x\ls 2\ve.\tag3.4$$

As (3.4) holds for any given $\ve\in(0,1)$, we must have $\Delta(x)=o(x)$ as desired. \qed

\medskip
\noindent{\it Proof of Theorem 1.3}.  For $z=-e^{2\pi i/m}$ we have $\Re(z)=-\cos(2\pi/m)<1$ as $m\not=2$. Combining (3.1) with (2.1), we obtain
$$\sum_{n\ls x}\mu_m(n)\l\lfloor\f xn\r\rfloor=F(1+z)x(\log x)^{z}+O\l(x(\log x)^{-1-\cos(2\pi/m)}\r).$$
By Lemma 3.2,
$$\sum_{n\ls x}\mu_m(x)\l\{\f xn\r\}=o(x).$$
Therefore
$$x\sum_{n\ls x}\f{\mu_m(n)}n=\sum_{n\ls x}\mu_m(n)\l(\l\lfloor\f xn\r\rfloor+\l\{\f xn\r\}\r)=F(1+z)x(\log x)^{z}+o(x)$$
and hence
$$\sum_{n\ls x}\f{\mu_m(n)}n=\Cal G(z)(\log x)^{z}+o(1)\tag3.5$$
since $F(1+z)=G(z)/z=\Cal G(z)$. Combining (3.5) with (1.10) and noting that $(\log x)^{z-1}\to0$ as $x\to\infty$, we get
$c_*(z)=0$. So (1.10) reduces to (1.15).

 For $m=5,6,\ldots$ we clearly have $\cos(2\pi/m)>0$ and hence (1.15) implies that $\sum_{n=1}^\infty\mu_m(n)/n=0$.
This concludes the proof. \qed

\Remark\ 3.1. The way we prove (1.14) can be modified to show the equality
$$\sum_{n=1}^\infty\f{\lambda(n)}n=0.\tag3.6$$
Since $\lambda=\nu_1^*$, we have $\sum_{n\ls x}\lambda(n)\{x/n\}=o(x)$ by Lemma 3.2.
So it suffices to prove
$\sum_{n\ls x}\lambda(n)\lfloor x/n\rfloor=o(x)$. In fact,
$$\align\sum_{d\ls x}\lambda(d)\l\lfloor\f xd\r\rfloor=&\sum_{d\ls x}\lambda(d)\sum_{q\ls x/d}1=\sum_{dq\ls x}\lambda(d)
=\sum_{n\ls x}\sum_{d\mid n}\lambda(d)
\\=&|\{1\ls n\ls x:\ n\ \t{is a square}\}|=\lfloor \sqrt x\rfloor=o(x).
\endalign$$

\heading{4. Proof of Theorem 1.4}\endheading

\proclaim{Lemma 4.1} Let $m\in\{1,5,6,\ldots\}$. Then the series
$$\zeta_m(1):=\sum_{n=1}^\infty\f{\nu_m(n)}n\quad\t{and}\quad\zeta_m^*(1):=\sum_{n=1}^\infty\f{\nu^*_m(n)}n$$
converge. Moreover, we have
$$\lim_{s\to1+}\zeta_m(s)=\zeta_m(1)\quad\t{and}\quad\lim_{s\to1+}\zeta_m^*(s)=\zeta_m^*(1).\tag4.1$$
\endproclaim
\Proof. Let $a(n)=\nu_m(n)$ for all $n\in\Z^+$, or $a(n)=\nu_m^*(n)$ for all $n\in\Z^+$. Set $A(x):=\sum_{n\ls x}a(n)$ for $x\gs1$,
and $f_s(t)=t^{-s}$ for $s\gs1$ and $t\gs2$. By the partial summation formula, for $x\gs x_0\gs2$ we have
$$\sum_{x_0<n\ls x}a(n)f_s(n)=A(x)f_s(x)-A(x_0)f_s(x_0)-\int_{x_0}^xA(t)f_s'(t)dt$$
and hence
$$\sum_{x_0<n\ls x}\f{a(n)}{n^s}=\f{A(x)}{x^s}-\f{A(x_0)}{x_0^s}+s\int_{x_0}^x\f{A(t)}{t^{s+1}}dt.\tag4.2$$
In view of (2.1) or (2.3) with $z=-e^{2\pi i/m}$, there is a constant $c>0$ depending on $m$ such that
$$|A(t)|\ls \f{ct}{(\log t)^{\cos(2\pi/m)+1}}\quad\t{for all}\ t\gs 2.\tag4.3$$
For any $s\gs1$, we have
$$\bg|\f{A(x)}{x^s}\bg|\ls\f{|A(x)|}x\ls \f{c}{(\log x)^{\cos(2\pi/m)+1}},\
\ \bg|\f{A(x_0)}{x_0^s}\bg|\ls\f{|A(x_0)|}{x_0}\ls \f{c}{(\log x_0)^{\cos(2\pi/m)+1}} $$
and
$$\align\bg|\int_{x_0}^x\f{A(t)}{t^{s+1}}dt\bg|\ls&\int_{x_0}^x\f{ct}{t^2}(\log t)^{-\cos(2\pi/m)-1}dt=\f{c(\log t)^{-\cos(2\pi/m)}}{-\cos(2\pi/m)}\bg|_{t=x_0}^x
\\=&\f{c}{\cos(2\pi/m)}\l(\f1{(\log x_0)^{\cos(2\pi /m)}}-\f1{(\log x)^{\cos(2\pi /m)}}\r)
\endalign$$
with the help of (4.3).

 Let $\ve>0$. Since $\cos(2\pi/m)>0$, by the above, there is an integer $N(\ve)\gs2$ such that if $x>x_0\gs N(\ve)$ then for any $s\gs1$ we have
$$\bg|\sum_{x_0<n\ls x}\f{a(n)}{n^s}\bg|\ls\bg|\f{A(x)}{x^s}\bg|+\bg|\f{A(x_0)}{x_0^s}\bg|+s\bg|\int_{x_0}^x\f{A(t)}{t^{s+1}}dt\bg|
\ls \f{\ve}2+\f{\ve}2+s\ve=(1+s)\ve.$$
Therefore the series $\sum_{n=1}^\infty a(n)/n^s$ converges for any $s\gs1$, in particular $\sum_{n=1}^\infty a(n)/n$ converges!

In view of the general properties of Dirichelt's series (cf. [T, p.\,291]), we immediately have
$$\lim_{s\to1+}\sum_{n=1}^\infty\f{a(n)}{n^s}=\sum_{n=1}^\infty\f{a(n)}n.$$
This concludes the proof. \qed
\medskip

\proclaim{Lemma 4.2} Let $m>4$ be an integer. For $\Re(s)>1$, we have
$$\f d{ds}\log\zeta_m(s)+e^{2\pi i/m}\f d{ds}\log\zeta(s)=v(s)\tag4.4$$
and
$$\f d{ds}\log\zeta_m^*(s)+e^{2\pi i/m}\f d{ds}\log\zeta(s)=v^*(s),\tag4.5$$
where $v(s)$ is a suitable holomorphic function in the region $\Re(s)>\log_2(2\cos\f{\pi}m)$
and $v^*(s)$ is a suitable holomorphic functions in the half plane $\Re(s)>1/2$.
\endproclaim
\Proof.  (i) Equation (4.5) can be proved in a way similar to the proof of [KY, Theorem 1]. Let $z=-e^{2\pi i/m}$ and
$$v^*(s)=\sum_p(\log p)\sum_{k=2}^\infty\f{z-z^k}{p^{ks}}\quad \t{for}\ \Re(s)>\f12.$$
If $\sigma=\Re(s)>1/2$, then
$|p^s|=p^{\sigma}\gs\sqrt2$ for any prime $p$, and hence
$$\align\sum_p(\log p)\bg|\sum_{k=2}^\infty\f{z-z^k}{p^{ks}}\bg|
\ls&\sum_p(\log p)\sum_{k=2}^\infty\f2{p^{k\sigma}}=\sum_p\f{2\log p}{p^{2\sigma}(1-p^{-\sigma})}
\\\ls&\f2{1-1/\sqrt2}\sum_p\f{\log p}{p^{2\sigma}}\ls\f{2\sqrt2}{\sqrt2-1}\sum_{n=1}^\infty\f{\log n}{n^{2\sigma}}<\infty.
\endalign$$
So $v^*(s)$ is a holomorphic function in the region $\Re(s)>1/2$.

When $\Re(s)>1$, in view of Euler's product $\prod_p(1-p^{-s})^{-1}=\zeta(s)$ and the formula (1.4), we have
$$\align &\f d{ds}\log\zeta_m^*(s)+e^{2\pi i/m}\f d{ds}\log\zeta(s)
\\=&-\sum_p\f d{ds}\log(1+e^{2\pi i/m}p^{-s})-e^{2\pi i/m}\sum_p\f d{ds}\log(1-p^{-s})
\\=&-\sum_p\f{e^{2\pi i/m}(-\log p)p^{-s}}{1+e^{2\pi i/m}p^s}-e^{2\pi i/m}\sum_p\f{-(-\log p)p^{-s}}{1-p^{-s}}
\\=&-\sum_p(\log p)\sum_{k=1}^\infty\l(\f{-e^{2\pi i/m}}{p^s}\r)^k-e^{2\pi i/m}\sum_p(\log p)\sum_{k=1}^\infty\f1{p^{sk}}
\\=&\sum_p(\log p)\sum_{k=2}^\infty\f{z-z^k}{p^{sk}}=v^*(s).
\endalign$$
This proves (4.5).

(ii) To prove (4.4), we set $z=-e^{2\pi i/m}$ and
$$v(s)=(z^2-z)\sum_p\f{\log p}{(p^s-1)(p^s-1+z)}\quad\t{for}\ \Re(s)>\log_2\l(2\cos\f{\pi}m\r).$$
Note that
$$|1-z|=\sqrt{\l(1+\cos\f{2\pi}m\r)^2+\l(\sin\f{2\pi}m\r)^2}=2\cos\f{\pi}m>2\cos\f{\pi}4=\sqrt2.$$
If $\sigma=\Re(s)>\log_2(2\cos\f{\pi}m)$, then for any prime $p$ we have
$|p^s|=p^{\sigma}\gs2^{\sigma}>2\cos\f{\pi}m=|1-z|$ and hence $p^s-1+z\not=0$.
For each prime $p>3$ and $\sigma=\Re(s)>\log_2(2\cos\f{\pi}m)>\f12$, as $p^{\sigma}>(2^{\sigma})^2>2$ we have
$$|p^s-1|\gs p^{\sigma}-1>\f{p^{\sigma}}2;$$
also,
$$|p^s-1+z|\gs p^{\sigma}-|1-z|>\l(1-\f1{\sqrt2}\r)p^{\sigma}$$
since $p^{\sigma}>(2^{\sigma})^2>|1-z|^2>|1-z|\sqrt2$.
As $\sum_p(\log p)/p^{2\sigma}$ converges for any $\sigma>1/2$, we see that
$$\sum_p\f{\log p}{(p^s-1)(p^s-1+z)}$$
converges absolutely in the half plane $\Re(s)>\log_2(2\cos\f{\pi}m)$. Therefore
$v(s)$ is indeed a holomorphic function in the region $\Re(s)>\log_2(2\cos\f{\pi}m)$.

When $\Re(s)>1$, using Euler's product $\prod_p(1-p^{-s})^{-1}=\zeta(s)$ and the formula (1.3) we get
$$\align &\f d{ds}\log\zeta_m(s)+e^{2\pi i/m}\f d{ds}\log\zeta(s)
\\=&\sum_p\f d{ds}\log\l(1-\f{e^{2\pi i/m}}{p^s-1}\r)-e^{2\pi i/m}\sum_p\f d{ds}\log(1-p^{-s})
\\=&\sum_p\f{(-z/(p^s-1)^2)p^s\log p}{1+z/(p^s-1)}+z\sum_p\f{-(-\log p)p^{-s}}{1-p^{-s}}
\\=&\sum_p(\log p)\(-\f{zp^s}{(p^s-1)(p^s-1+z)}+\f z{p^s-1}\)
\\=&\sum_p\f{(z^2-z)\log p}{(p^s-1)(p^s-1+z))}=v(s).
\endalign$$
This proves (4.4).

In view of the above, we have completed the proof of Lemma 4.2. \qed

\medskip

\noindent{\it Proof of Theorem 1.4}.
 Let $m\in\{1,3,4,\ldots\}$ and $z=-e^{2\pi i/m}$. When $m=3$, (1.9) and (1.11) yield (1.16) with $c_3=c(z)$
 and $C_3=C(z)$. In the case $m=4$, (1.9) and (1.11) give (1.17) with $c_4=c(-i)$ and $C_4=C(-i)$.

Now we assume that $m=1$ or $m>4$. Note that $\Re(z)=-\cos(2\pi/m)<0$. By (1.9) and (1.11), we have
$$V_m(x)=\Cal F(z)(\log x)^z+c_m+O((\log x)^{z-1})$$
and $$V^*_m(x)=\Cal H(z)(\log x)^z+C_m+O((\log x)^{z-1}),$$
where $c_m=c(z)$ and $C_m=C(z)$. It follows that
$$\lim_{x\to\infty}V_m(x)=c_m\quad\t{and}\quad\lim_{x\to\infty}V_m^*(x)=C_m.$$
Also, (1.19) and (1.20) hold if $c_m=C_m=0$. So it suffices to show
$V_m(x)=o(1)$ and $V_m^*(x)=o(1)$. This holds for $m=1$ since $\zeta_1^*(1)=\sum_{n=1}^\infty\lambda(n)/n=0$ by (3.6)
and $\zeta_1(1)=\sum_{n=1}^\infty(-1)^{\omega(n)}/n=0$ by [LD].

Below we fix $m\in\{5,6,\ldots\}$. By Lemma 4.2, we have
$$\f d{ds}\log\l(\zeta_m(s)\zeta(s)^{e^{2\pi i/m}}\r)=v(s)\quad\t{for}\ s>1,$$
where $v(s)$ is a holomorphic function in the half plane $\Re(s)>\log_2(2\cos\f{\pi}m)$.
Choose a number $s_0>1$. Then
$$\zeta_m(s)\zeta(s)^{e^{2\pi i/m}}=\zeta_m(s_0)\zeta(s_0)^{e^{2\pi i/m}}e^{\int_{s_0}^s v(t)dt}$$
for all $s>1$, and hence
$$\lim_{s\to1+}\zeta_m(s)=\lim_{s\to 1+}\l(\f{\zeta(s_0)}{\zeta(s)}\r)^{e^{2\pi i/m}}\zeta_m(s_0)e^{\int_{s_0}^1 v(t)dt}=0$$
since $1>\log_2(2\cos\f{\pi}m)$, $\Re(e^{2\pi i/m})=\cos\f{2\pi}m>0$ and $\lim_{s\to1+}\zeta(s)=\infty$.
Similarly, by applying (4.5) in Lemma 4.2 we get $\lim_{s\to1+}\zeta_m^*(s)=0$.
Combining these with Lemma 4.1 we finally obtain
$$\zeta_m(1)=\lim_{s\to1+}\zeta_m(s)=0\ \ \ \t{and}\ \ \ \zeta_m^*(1)=\lim_{s\to1+}\zeta_m^*(s)=0$$
as desired. This concludes the proof of Theorem 1.4. \qed

\heading{5. Proof of Theorem 1.5}\endheading

\medskip
\noindent{\it Proof of Theorem 1.5}.
Let $L(x)=\sum_{n\ls x}(-1)^{\Omega(n)}$. Formula (2.3) with $z=-1$ yields that $L(x)=o(x)$. Observe that
$$S(x)+L(x)=\sum_{n\ls x}((-1)^n+1)(-1)^{\Omega(n)}=2\sum_{m\ls x/2}(-1)^{\Omega(2m)}=-2L\l(\f x2\r).$$
Therefore
$$S(x)=-L(x)-2L\l(\f x2\r)=o(x).$$

For any complex number $s$, obviously
$$\sum_{n\ls x}\f{(-1)^{n-\Omega(n)}}{n^s}+\sum_{n\ls x}\f{\lambda(n)}{n^s}=2\sum\Sb n\ls x\\2\mid n\endSb\f{\lambda(n)}{n^s}
=-2\sum_{m\ls x/2}\f{\lambda(m)}{(2m)^s}$$
and hence
$$\sum_{n\ls x}\f{(-1)^{n-\Omega(n)}}{n^s}=-2^{1-s}\sum_{n\ls x/2}\f{\lambda(n)}{n^s}-\sum_{n\ls x}\f{\lambda(n)}{n^s}.$$
Since $\sum_{n\ls x}\lambda(n)/n=o(1)$ by (3.6), we get $\sum_{n\ls x}(-1)^{n-\Omega(n)}/n=o(1)$
and hence $\sum_{n=1}^\infty(-1)^{n-\Omega(n)}/n=0$.

Let $\Re(s)>1$. Note that
$$\sum_{n=1}^\infty\f{(-1)^{n-\Omega(n)}}{n^s}=-(1+2^{1-s})\sum_{n=1}^\infty\f{\lambda(n)}{n^s}=-(1+2^{1-s})\f{\zeta(2s)}{\zeta(s)}.$$
On the other hand, by the partial summation method, we have
$$\sum_{n\ls x}\f{(-1)^{n-\Omega(n)}}{n^s}=\f{S(x)}{x^s}+s\int_1^x\f{S(t)}{t^{s+1}}dt$$
and hence
$$\sum_{n=1}^\infty\f{(-1)^{n-\Omega(n)}}{n^s}=s\int_1^\infty \f{S(t)}{t^{s+1}}dt.$$
Therefore
$$-(1+2^{1-s})\f{\zeta(2s)}{\zeta(s)}=s\int_1^\infty \f{S(t)}{t^{s+1}}dt.\tag5.1$$

Let $\sigma_c$ be the least real number such that the integral in (5.1) converges whenever $\Re(s)>\sigma_c$.
By the above, $\sigma_c\ls 1$.

Suppose that $S(x)>0$ for all $x\gs 5$. In view of (5.1), by applying Landau's theorem (cf. [MV, Lemma 15.1] or Ex. 16 of [Ap, p.248])
 we obtain
$$\lim_{s\to\sigma_c+}-\f{1+2^{1-s}}s\cdot\f{\zeta(2s)}{\zeta(s)}=\infty$$
and hence $\sigma_c\ls1/2$ since $\zeta(s)$ has no real zeroes with $s>1/2$. (Note that $(1-2^{1-s})\zeta(s) =\sum_{n=1}^\infty(-1)^{n-1}/n^s\not=0$
for all $s>0$ with $s\not=1$.) So the right-hand side of (5.1) converges for $\Re(s)>1/2$
and hence so is the left-hand side of (5.1). Therefore $\zeta(s)\not=0$ for $\Re(s)>1/2$, i.e., the Riemann Hypothesis holds.

Similarly, if $T(x)<0$ for all $x\gs1$, then we get the Riemann Hypothesis by applying Landau's theorem.

So far we have completed the proof of Theorem 1.5. \qed
\medskip

\Ack. The author would like to thank Dr. D. Broadhurst, P. Humphries, S. Kim, W. Narkiewicz, H. Pan, M. Radziwill, P. Xi, L.-L. Zhao
and the referee for helpful comments.


\widestnumber\key{BFM}

 \Refs

\ref\key A\by A. W. Addison\paper A note on the compositeness of numbers\jour Proc. Amer. Math. Soc.
\vol \yr 1957\pages 151--154\endref


\ref\key Ap\by T. M. Apostol\book Introduction to Analytic Number Theory\publ Sptinger, New York, 1976\endref

\ref\key BFM\by P. Borwein, R. Ferguson and M. J. Mossinghoff\paper Sign changes in sums of the Liouville function
\jour Math. Comp.\vol 77\yr 2008\pages 1681-1694\endref

\ref\key BS\by R. Br\"oker and A. V. Sutherland\paper An explicit height bound for the classical modular polynomial\jour Ramanujan J.\vol 22\yr 2010\pages 293--313\endref

\ref\key CD\by M. Coons and S. R. Dahmen\paper On the residue class distribution of the number of prime divisors of an integer
\jour Nagoya Math. J. \vol 202\yr 2011\pages 15--22\endref

\ref\key CP\by R. Crandall and C. Pomerance\book Prime Numbers: A Computational Perspective\publ 2nd Edition, Springer, New York, 2005\endref

\ref\key D\by H. Delange\paper Sur des formules de Atle Selberg\jour Acta Arith.\vol 19\yr 1971\pages 105-146\endref

\ref\key G\by E. Grosswald\paper The average order of an arithmetic function
\jour Duke Math. J.\vol 23\yr 1956\pages 41--44\endref

\ref\key H\by C. B. Haselgrove\paper A disproof of a conjecture of P\'olya\jour Mathematika\vol 5\yr 1958\pages 141--145\endref

\ref\key IR\by K. Ireland and M. Rosen\book A Classical Introduction to Modern Number Theory, 2nd Edition
\publ Springer, New York, 1990\endref

\ref\key KY\by T. Kudota and M. Yoshida\paper A note on the congruent distribution of the number of prime factors of natural numbers
\jour Nagoya Math. J.\vol 163\yr 2001\pages 1--11\endref

\ref\key L\by R. S. Lehman\paper On Liouville's function\jour Math. Comp.\vol 14\yr 1960\pages 311--320\endref

\ref\key LD\by J. van de Lune and R. E. Dressler\paper Some theorems concerning the number theoretic function $\omega(n)$
\jour J. Reine Angew. Math.\vol 277\yr 1975\pages 117--119\endref

\ref\key MV\by H. L. Montgomery and R. C. Vaughan\book Multiplicative Number Theory I. Classical Theory
\publ Cambridge Univ. Press, Cambridge, 2007\endref

\ref\key S\by A. Selberg\paper Note on a paper by L. G. Sathe\jour J. Indian Math. Soc. \vol 18\yr 1954\pages 83-87\endref

\ref\key T\by E. C. Titchmarsh\book The Theory of Functions\publ  Reprint of the second (1939) edition, Oxford University Press, Oxford, 1958\endref
\endRefs
\enddocument